\documentclass[12pt]{article}

\newcommand{\be}{\begin{equation}}
\newcommand{\ee}{\end{equation}}

\usepackage{epsfig,amsmath} 
\usepackage{graphics}

\begin{document}

\begin{center}
Journal of Physics A. Vol.38. No.26. (2005) pp.5929-5943. \\
\end{center}

\begin{center}
{\Large \bf Fractional Generalization of Gradient and \\
Hamiltonian Systems}
\end{center}

\vskip 7mm
\begin{center}
{\large \bf Vasily E. Tarasov } \\
{\it Skobeltsyn Institute of Nuclear Physics, \\ 
Moscow State University, Moscow 119992, Russia \\
E-mail: tarasov@theory.sinp.msu.ru } \\
\end{center}

\begin{abstract}
We consider a fractional generalization of 
Hamiltonian and gradient systems. 
We use differential forms 
and exterior derivatives of fractional orders. 
We derive fractional generalization of Helmholtz conditions
for phase space.
Examples of fractional gradient and 
Hamiltonian systems are considered. 
The stationary states for these systems are derived. 
\end{abstract}

\vskip 3 mm
PACS numbers: 45.20.-d; 05.45.-a;
\vskip 5mm

\section{Introduction}

Derivatives and integrals of fractional order \cite{SKM,OS} have 
found many applications in recent studies in physics.
The interest in fractional analysis has been growing continually 
during the past few years. 
Fractional analysis has numerous applications: 
kinetic theories \cite{Zaslavsky1,Zaslavsky2,Physica},  
statistical mechanics \cite{chaos,PRE05,JPCS}, 
dynamics in complex media \cite{Nig,PLA05,PLA05-2,AP05,Chaos05}, 
and many others \cite{M-K1,M-K2,Hilfer,Carpinteri-Mainardi}.

The theory of derivatives of non-integer order goes back 
to Leibniz, Liouville, Grunwald, Letnikov and Riemann. 
In the past few decades many authors have pointed out that 
fractional-order models are more appropriate than integer-order 
models for various real materials. 
Fractional derivatives provide an excellent instrument 
for the description of memory and hereditary properties 
of various materials and processes. 
This is the main advantage of fractional derivatives 
in comparison with classical integer-order models in 
which such effects are, in fact, neglected. 
The advantages of fractional derivatives become apparent 
in modelling mechanical and electrical properties of real materials, 
as well as in the description of rheological properties of rocks, 
and in many other fields. 

In this paper, we use
a fractional generalization of exterior calculus that was suggested
in \cite{FDF,FDF2}. Fractional generalizations of differential forms
and exterior derivatives were defined in \cite{FDF}.
It allows us to consider the fractional generalization
of Hamiltonian and gradient dynamical systems \cite{Gilmor,DNF}.
The suggested class of fractional gradient and Hamiltonian systems
is wider than the usual class of gradient and Hamiltonian 
dynamical systems. 
The gradient and Hamiltonian systems can be considered as
a special case of fractional gradient and Hamiltonian systems. 

In section 2, a brief review of 
gradient systems and exterior calculus 
is considered to fix notation and provide a convenient reference.
In section 3, a brief review of fractional (exterior) calculus
will be given to fix notations and provide a convenient reference.
In section 4, a definition of fractional generalization 
of gradient systems is suggested.
In section 5, we consider a fractional gradient system that
cannot be considered as a gradient system. 
In section 6, we prove that a dynamical system that 
is defined by the well-known Lorenz equations \cite{Lor1,Lor2} 
can be considered as a fractional gradient system.
In section 7, a brief review of Hamiltonian systems  
is considered to fix notations and provide a convenient reference.
In section 8, we consider the fractional generalization
of Hamiltonian systems and Helmholtz conditions. 
In section 9, the simple example of fractional Hamiltonian systems 
is discussed.
Finally, a short conclusion is given in section 10.

\section{Gradient systems}

In this section, a brief review of 
gradient systems and exterior calculus \cite{DNF}
is considered to fix notations and provide a convenient reference.

Gradient systems arise in dynamical systems theory \cite{Gilmor,HS,DNF}.
They are described by the equation
$d{\bf x}/dt=-grad V(x)$, where ${\bf x} \in R^n$.
In Cartesian coordinates, the gradient is given by
$grad V={\bf e}_i {\partial V}/{\partial x_i}$,
where ${\bf x}={\bf e}_i x_i$. 
Here and later, we mean the sum on the repeated indices 
$i$ and $j$ from 1 to n.  \\

{\bf Definition 1.}
{\it A dynamical system that is described by the equations
\be  \label{cs} \frac{dx_i}{dt}=F_i(x)  \quad (i=1,...,n) \ee
is called a gradient system in $R^n$ if the differential 1-form 
\be \label{omega} \omega=F_i(x)dx_i \ee
is an exact form $\omega=-dV$, where
$V=V(x)$ is a continuously differentiable function (0-form). } \\

Here $d$ is the exterior derivative \cite{DNF}.
Let $V=V(x)$ be a real, continuously differentiable function on $R^n$.
The exterior derivative of the function $V$ is the 1-form
$dV=dx_i {\partial V}/{\partial x_i}$ 
written in a coordinate chart $(x_1,...,x_n)$.

In mathematics \cite{DNF}, 
the concepts of closed form and exact form are defined 
for differential forms by the equation $d\omega =0$
for a given form $\omega$ to be a closed form 
and $\omega=dh$ for an exact form.
It is known, that to be exact is a sufficient condition to be closed. 
In abstract terms, the question of whether this is also a necessary condition 
is a way of detecting topological information, by differential conditions.

Let us consider the 1-form (\ref{omega}). 
The formula for the exterior derivative $d$ of differential form 
(\ref{omega}) is
\[ d \omega=\frac{1}{2} \left(\frac{\partial F_i}{\partial x_j} - 
\frac{\partial F_j}{\partial x_i} \right) dx_j \wedge dx_i , \]
where $\wedge$ is the wedge product. 
Therefore, the condition for $\omega$ to be closed is
\[ \frac{\partial F_i}{\partial x_j} - 
\frac{\partial F_j}{\partial x_i}=0. \]
In this case, if $V(x)$ is a potential function then
$dV = dx_i{\partial V}/{\partial x_i}$.
The implication from 'exact' to 'closed' is then a consequence 
of the symmetry of the second derivatives: 
\be \label{Vxy} \frac{\partial^2 V}{\partial x_i \partial x_j}=
\frac{\partial^2 V}{\partial x_j \partial x_i} . \ee
If the function $V=V(x)$ is smooth function, then the second derivative
commute, and equation (\ref{Vxy}) holds. 

The fundamental topological result here is the Poincare lemma. 
It states that for a contractible open subset $X$ of $R^n$, 
any smooth p-form $\beta$ defined on $X$ that is closed, 
is also exact, for any integer $p > 0$ 
(this has content only when $p$ is at most $n$).
This is not true for an open annulus in the plane, for some 
1-forms $\omega$ that fail to extend smoothly to the whole disk, 
so that some topological condition is necessary.
A space $X$ is contractible if the identity map on $X$ is homotopic 
to a constant map. Every contractible space is simply connected.
A space is simply connected if it is path connected and every loop 
is homotopic to a constant map. \\

{\bf Proposition 1.}
{\it If a smooth vector field ${\bf F}={\bf e}_i F_i(x)$ of system (\ref{cs})
satisfies the relations 
\be \label{Potcond}
\frac{\partial F_i}{\partial x_j}-\frac{\partial F_j}{\partial x_i}=0 \ee
on a contractible open subset X of $R^n$,
then the dynamical system (\ref{cs}) is
the gradient system such that }
\be \label{ps} \frac{dx_i}{dt}=-\frac{\partial V(x)}{\partial x_i} . \ee

This proposition is a corollary of the Poincare lemma.
The Poincare lemma states that for a contractible open 
subset $X$ of $R^n$, any smooth 1-form (\ref{omega})
defined on X that is closed, is also exact.


The equations of motion for the gradient system 
on a contractible open subset X of $R^n$ 
can be represented in the form (\ref{ps}).
Therefore, the gradient systems can be defined by 
the potential function $V=V(x)$.

If the exact differential 1-form $\omega$ is equal to zero
($dV=0$), then we get the equation 
\be \label{VC} V(x)-C=0, \ee
which defines the stationary states of the gradient dynamical 
system (\ref{ps}). Here $C$ is a constant.

\section{Fractional differential forms}

If the partial derivatives in the definition of the exterior derivative 
$d=dx_i \partial / \partial x_i$ are allowed to assume fractional order, 
a fractional exterior derivative can be defined \cite{FDF} by the equation
\be \label{fed}  d^{\alpha}=(dx_i)^{\alpha} {\bf D}^{\alpha}_{x_i} . \ee
Here, we use the fractional derivative ${\bf D}^{\alpha}_{x}$
in the Riemann-Liouville form \cite{SKM} that is defined by the equation
\be \label{df} {\bf D}^{\alpha}_{x}f(x)=\frac{1}{\Gamma (m-\alpha)}
\frac{\partial^m}{\partial x^m}  
\int^x_{0} \frac{f(y) dy}{(x-y)^{\alpha-m+1}} , \ee
where $m$ is the first whole number greater than  or equal to $\alpha$.
The initial point of the fractional derivative \cite{SKM} is set to zero.
The derivative of powers $k$ of $x$ is
\be \label{xk} {\bf D}^{\alpha}_x x^k=
\frac{\Gamma(k+1)}{\Gamma(k+1-\alpha)} x^{k-\alpha} ,\ee
where $k \ge 1$, and $\alpha \ge 0$. 
The derivative of a constant $C$ need not be zero
\be \label{const}
{\bf D}^{\alpha}_x C=\frac{x^{-\alpha}}{\Gamma(1-\alpha)} C. \ee

For example, the fractional exterior derivative 
of order $\alpha$ of $x^k_1$, with the initial point taken 
to be zero and $n=2$, is given by
\be \label{daq} d^{\alpha} x^k_1= (dx_1)^{\alpha} {\bf D}^{\alpha}_{x_1} x^k_1+
(dx_2)^{\alpha} {\bf D}^{\alpha}_{x_2} x^k_1 . \ee
Using  equation (\ref{xk}), we get the following relation 
for the fractional exterior derivative of $x^k_1$: 
\[ d^{\alpha} x^k_1= (dx_1)^{\alpha} 
\frac{\Gamma(k+1) x^{k-\alpha}_1}{\Gamma(k+1-\alpha)} +
(dx_2)^{\alpha} \frac{x^{k}_1 x^{-\alpha}_2}{\Gamma(1-\alpha)} . \]

\section{Fractional gradient systems}

A fractional generalization of exterior calculus was suggested
in \cite{FDF,FDF2}. A fractional exterior derivative 
and the fractional differential forms were defined \cite{FDF}.
It allows us to consider the fractional generalization
of gradient systems.
 
Let us consider a dynamical system  that is defined
by the equation ${d{\bf x}}/{dt}={\bf F}$, 
on a subset X of $R^n$. 
In Cartesian coordinates, we can use the following equation
\be \label{cs2}  \frac{dx_i}{dt}=F_i(x), \ee
where $i=1,..,n$, ${\bf x}={\bf e}_i x_i$, and ${\bf F}={\bf e}_i F_i(x)$.
The fractional analog of definition 1 has the form. \\

{\bf Definition 2.}
{\it A dynamical system (\ref{cs2})
is called a fractional gradient system if the 
fractional differential 1-form 
\be \label{omega-a} \omega_{\alpha}=F_i(x)(dx_i)^{\alpha} \ee 
is an exact fractional form $\omega_{\alpha}=-d^{\alpha} V$, 
where $V=V(x)$ is a continuously differentiable function. } \\

Using the definition of the fractional exterior derivative, 
equation (\ref{omega-a}) can be represented as
\[ \omega_{\alpha}=-d^{\alpha} V=-(dx_i)^{\alpha} {\bf D}^{\alpha}_{x_i}V. \]
Therefore, we have $F_i(x)=-{\bf D}^{\alpha}_{x_i} V$.

Note that equation (\ref{omega-a}) is a fractional generalization
of equation (\ref{omega}). 
If $\alpha=1$, then equation (\ref{omega-a}) leads us to
equation (\ref{omega}). 
Obviously, a fractional 1-form $\omega_{\alpha}$
can be closed when the 1-form $\omega=\omega_1$ is not closed. 
The fractional generalization of the Poincare lemma is considered
in \cite{FDF2}. Therefore, we have the following proposition. \\

{\bf Proposition 2.}
{\it If a smooth vector field ${\bf F}={\bf e}_i F_i(x)$ 
on a contractible open subset X of $R^n$ 
satisfies the relations
\be \label{FHC} {\bf D}^{\alpha}_{x_j} F_i- {\bf D}^{\alpha}_{x_i} F_j=0 , \ee
then the dynamical system (\ref{cs2})
is a fractional gradient system such that 
\be \label{fps} \frac{dx_i}{dt}=-{\bf D}^{\alpha}_{x_i} V(x) , \ee
where $V(x)$ is a continuous differentiable function and 
${\bf D}^{\alpha}_{x_i}V=-F_i$. }\\

{\bf Proof.} This proposition is a corollary of the 
fractional generalization of Poincare lemma \cite{FDF2}.
The Poincare lemma is shown \cite{FDF,FDF2}
to be true for the exterior fractional derivative. \\

Relations (\ref{FHC}) are the fractional generalization
of relations (\ref{Potcond}). 
Note that the fractional derivative of a 
constant need not be zero (\ref{const}). 
Therefore, we see that constants $C$ in the equation
$V(x)=C$ cannot define a stationary state of the 
gradient system (\ref{fps}). 
It is easy to see that 
\[ {\bf D}^{\alpha}_{x_i} V(x)={\bf D}^{\alpha}_{x_i} C=
\frac{x^{-\alpha}_i}{\Gamma(1-\alpha)} C \not=0. \]
In order to define stationary states of fractional gradient systems, 
we consider the solutions of the system of the equations 
${\bf D}^{\alpha}_{x_i} V(x)=0$. \\

{\bf Proposition 3.}
{\it The stationary states of gradient system (\ref{fps})
are defined by the equation 
\be \label{ssgs} V(x)-|\prod^n_{i=1} x_i|^{\alpha-m} 
\sum^{m-1}_{k_1=0} ... \sum^{m-1}_{k_n=0} C_{k_1...k_n} 
\Bigl( \prod^n_{i=1}(x_i)^{k_i} \Bigr)=0 . \ee
The $C_{k_1...k_n}$ are constants and $m$ is the first whole number
greater than  or equal to $\alpha$.}\\

{\bf Proof.} 
In order to define the stationary states of a fractional gradient system, 
we consider the solution of the equation
\be \label{DH} {\bf D}^{\alpha}_{x_i} V(x)=0 . \ee
This equation can be solved by using equation (\ref{df}).
Let $m$ be the first whole number greater than or equal to $\alpha$;
then we have the solution \cite{SKM,OS} of equation (\ref{DH}) in the form
\be \label{sDF} V(x)=|x_i|^{\alpha} \sum^{m-1}_{k=0} 
a_k(x_1,...,x_{i-1},x_{i+1},...,x_n) (x_i)^{k} , \ee
where $a_k$ are functions of the other coordinates.
Using equation (\ref{sDF}) for $i=1,...,n$,  
we get the solution of the system of equation (\ref{DH})
in the form (\ref{ssgs}). \\

If we consider $n=2$ such that $x=x_1$ and $y=x_2$, we have
the equations of motion for fractional gradient system
\be \label{xH2} \frac{dx}{dt} =-{\bf D}^{\alpha}_{x} V(x,y), 
\quad \frac{dy}{dt} =-{\bf D}^{\alpha}_{y} V(x,y) .\ee
The stationary states of this system are defined by the equation 
\[ V(x,y)-|xy|^{\alpha-1} \sum^{m-1}_{k=0} \sum^{m-1}_{l=0} C_{kl} 
x^{k} y^{l}=0 . \]
The $C_{kl}$ are constants and $m$ is the first whole number
greater than  or equal to $\alpha$.

\section{Examples of fractional gradient system}

In this section, we consider a fractional gradient systems that
cannot be considered as a gradient system. 
We prove that the class of fractional gradient systems
is wider than the usual class of gradient dynamical systems. 
The gradient systems can be considered as
a special case of fractional gradient systems. \\

{\bf Example 1.}
Let us consider the dynamical system that is defined by the equations
\be \label{eqex1} \frac{dx}{dt}=F_x, \quad \frac{dy}{dt}=F_y ,\ee
where the right hand sides have the form
\be \label{ex1} 
F_x=a c x^{1-k}+ b x^{-k} , \quad F_y=(ax+b) y^{-k} , \ee
where $a\not=0$. 
This system cannot be considered as a gradient dynamical system.
Using 
\[ \frac{\partial F_x}{\partial y}-\frac{\partial F_y}{\partial x}=
ay^{-k} \not=0 , \]
we get that $\omega=F_xdx+F_ydy$ is not closed form
\[ d\omega=-a y^{-k} dx \wedge dy .\]
Note that the relation (\ref{FHC}) in the form
\[ {\bf D}^{\alpha}_y F_x-{\bf D}^{\alpha}_x F_y=0 , \]
is satisfied for the system (\ref{ex1}), if $\alpha=k$ and 
the constant $c$ is defined by
\[ c=\frac{\Gamma(1-\alpha)}{\Gamma(2-\alpha)} . \]
Therefore, this system
can be considered as a fractional gradient system with 
the linear potential function
\[ V(x,y)=\Gamma(1-\alpha) (ax+b) ,\]
where $\alpha=k$. \\

{\bf Example 2.}
Let us consider the dynamical system that is defined 
by equation (\ref{eqex1}) with  
\be F_x=an(n-1)x^{n-2}+ck(k-1)x^{k-2}y^l , \ee
\be F_y=bm(m-1)y^{m-2}+cl(l-1)x^ky^{l-2} , \ee
where $k\not=1$ and $l\not=1$. It is easy to derive that 
\[ \frac{\partial F_x}{\partial y}-\frac{\partial F_y}{\partial x}=
ckl \ x^{k-2} y^{l-2} [(k-1)y-(l-1)x] \not=0 ,  \]
and the differential form $\omega=F_xdx+F_ydy$ is not closed $d \omega\not=0$.
Therefore, this system is not a gradient dynamical system. 
Using conditions (\ref{FHC}) in the form
\[ {\bf D}^2_y F_x-{\bf D}^2_x F_y=
\frac{\partial^2 F_x}{\partial y^2}-
\frac{\partial^2 F_y}{\partial x^2}=0 ,  \]
we get $d^{\alpha} \omega=0$ for $\alpha=2$. 
As the result, we have that this system can be considered as 
a fractional gradient system with the potential function
\[ V(x,y)=ax^n+by^m+cx^ky^l . \]

In the general case, the fractional gradient system cannot
be considered as a gradient system. 
The gradient systems can be considered as a special case 
of fractional gradient systems such that $\alpha=1$.

\section{Lorenz system as a fractional gradient system}

In this section, we prove that dynamical systems that 
are defined by the well-known Lorenz equations \cite{Lor1,Lor2} 
are fractional gradient system. \\

The well-known Lorenz equations \cite{Lor1,Lor2} 
are defined by 
\[ \frac{dx}{dt}=F_x ,\quad 
\frac{dy}{dt}=F_y , \quad 
\frac{dz}{dt}=F_z , \]
where the right hand sides $F_x$, $F_y$ and $F_z$ have the forms
\[ F_x=\sigma (y-x) ,\quad F_y=(r-z)x-y , \quad  F_z=xy-bz . \]
The parameters $\sigma$, $r$ and $b$ can be equal to the following 
values
\[ \sigma=10, \quad b=8/3, \quad r=470/19 \simeq 24.74 \ . \]
The dynamical system which is defined by the Lorenz equations
cannot be considered as a gradient dynamical system. 
It is easy to see that
\[ \frac{\partial F_x}{\partial y}-\frac{\partial F_y}{\partial x}
=z+\sigma-r ,\]
\[ \frac{\partial F_x}{\partial z}-\frac{\partial F_z}{\partial x}
=-y , \quad \quad \quad \quad \]
\[ \frac{\partial F_y}{\partial z}-\frac{\partial F_z}{\partial y}
=-2x . \quad \quad\quad \]
Therefore, $\omega=F_xdx+F_ydy+F_zdz$ is not a closed 1-form and we have 
\[ d \omega=-(z+\sigma-r)dx \wedge dy+y dx \wedge dz+
2x dy \wedge dz . \]

For the Lorenz equations, conditions (\ref{FHC}) 
can be satisfied in the form
\[ {\bf D}^2_y F_x-{\bf D}^2_x F_y=0,
\quad {\bf D}^2_z F_x-{\bf D}^2_x F_z=0,
\quad {\bf D}^2_z F_y-{\bf D}^2_y F_z=0 . \]
As the result, we get that the Lorenz system can be considered 
as a fractional gradient dynamical system with 
potential function
\be \label{LP} V(x,y,z)=\frac{1}{6}\sigma x^3-\frac{1}{2}\sigma yx^2
+\frac{1}{2}(z-r)xy^2+ \frac{1}{6} y^3 
-\frac{1}{2}xyz^2+\frac{b}{6}z^3 . \ee
The potential (\ref{LP}) uniquely defines the Lorenz system.
Using equation (\ref{ssgs}), we can get that the stationary states of 
the Lorenz system are defined by the equation
\be \label{LSS} V(x,y,z)+
C_{00}+C_x x+C_y y+C_z z+C_{xy}xy+C_{xz}xz+C_{yz}yz=0, \ee
where $C_{00}$, $C_x$, $C_y$, $C_z$ $C_{xy}$, $C_{xz}$, and $C_{yz}$
are the constants and $\alpha=m=2$. 
The plot of stationary states of Lorenz system
with the following constants 
$C_{00}=1$, $C_x=C_y=C_z=C_{xy}=C_{xz}=C_{yz}=0$ and
parameters $\sigma=10$, $b=3$, and $r=25$ is shown in figure 1 and 2.

Note that 
the Rossler system \cite{Ros}, which is defined by the equations
\[ \frac{dx}{dt}=-(y+z), \quad 
\frac{dx}{dt}=x+0.2y, \quad \frac{dz}{dt}=0.2+(x-c)z ,\] 
can be considered as a fractional gradient system with 
the potential function
\be \label{RPF} V(x,y,z)=\frac{1}{2}(y+z)x^2-\frac{1}{2}xy^2-
\frac{1}{30} y^3 -\frac{1}{10}z^2-\frac{1}{6}(x-c)z^3 . \ee
This potential uniquely defines the Rossler system. 
The stationary states of the Rossler system 
are defined by equation (\ref{LSS}), where the potential function
is defined by (\ref{RPF}). 
The plot of stationary states of the Rossler system
for the constants
$C_{00}=1$, $C_x=C_y=C_z=C_{xy}=C_{xz}=C_{yz}=0$ and
parameter $c=1$ is shown in figures 3 and 4.

Let us note the interesting qualitative property of surfaces (\ref{LSS}) 
which is difficult to see from figures.
The surfaces of the stationary states of the Lorenz and Rossler systems
separate the three dimensional Euclidean space into some number of areas. 
We have eight areas for the Lorenz system, and
four areas for the Rossler system.
This separation has the interesting property for
some values of parameters.
All regions are connected with each other. 
Beginning movement from one of the areas, it is possible 
to appear in any other area, not crossing a surface.
Any two points from different areas can be connected 
by a curve which does not cross a surface.
It is difficult to see this property from the figures 1-4.

\section{Hamiltonian systems}

In this section, a brief review of Hamiltonian systems  
is considered to fix notations and provide a convenient reference.

Let us consider the canonical coordinates 
$(q^1,...,q^n,p^{1},...,p^{n})$
in the phase space $R^{2n}$. 
We consider a dynamical system that is defined by the equations
\be \label{eq1} \frac{dq_i}{dt}=G^i(q,p), 
\quad \frac{dp_i}{dt}=F^i(q,p) . \ee
The definition of Hamiltonian systems can be realized in the 
following form \cite{Tartmf3,JPA05}. \\

{\bf Definition 3.} 
{\it A dynamical system (\ref{eq1}) on the phase space $R^{2n}$, 
is called a Hamiltonian system if the differential 1-form 
\be \label{beta} \beta=G^idp_i-F^idq_i, \ee
is a closed form $d \beta=0$, where $d$ is the exterior derivative. 
A dynamical system is called a non-Hamiltonian system 
if the differential 1-form $\beta$ is nonclosed $d \beta \not=0$.} \\

The exterior derivative for the phase space is defined as
\be \label{ed}  d= dq_i \frac{\partial}{\partial q_i}+
dp_i \frac{\partial}{\partial p_i} . \ee
Here and later, we mean the sum on the repeated indices
$i$ and $j$ from 1 to n. \\

{\bf Proposition 4.} 
{\it If the right-hand sides of equations (\ref{eq1})
satisfy the Helmholtz conditions \cite{Helm,Tartmf3,JPA05} for the 
phase space, which have the following forms: 
\be \label{HC1} \frac{\partial G^{i}}{\partial p_j}-
\frac{\partial G^{j}}{\partial p_i}= 0, \ee
\be \label{HC2} \frac{\partial G^{j}}{\partial q_i}+
\frac{\partial F^{i}}{\partial p_j}=0, \ee
\be \label{HC3} \frac{\partial F^{i}}{\partial q_j}-
\frac{\partial F^{j}}{\partial q_i}= 0, \ee
then the dynamical system (\ref{eq1}) is a Hamiltonian system.} \\

{\bf Proof}. 
In the canonical coordinates $(q,p)$, the vector fields that 
define the system have the components $(G^i,F^i)$, 
which are used in equation (\ref{eq1}). 
Let us consider the 1-form that is defined by the equation
\[ \beta=G^idp_i-F^idq_i. \]
The exterior derivative for this form can be written by the relation:
\[ d\beta=d(G^idp_i)-d(F^idq_i).  \]
It now follows that
\[ d\beta= \frac{\partial G^i}{\partial q_j} 
dq_j \wedge dp_i
+\frac{\partial G^i}{\partial p_j} dp_j \wedge dp_i-
\frac{\partial F^i}{\partial q_j} dq_j \wedge dq_i-
\frac{\partial F^i}{\partial p_j} dp_j \wedge dq_i  . \]
Here, $\wedge$ is the wedge product. 
This equation can be rewritten in an equivalent form as
\[ d\beta=
\left( \frac{\partial G^j}{\partial q_i} 
+\frac{\partial F^i}{\partial p_j}\right) dq_i \wedge dp_j+
\frac{1}{2}\left( \frac{\partial G^j}{\partial p_i}
-\frac{\partial G^i}{\partial p_j} \right)dp_i \wedge dp_j
+\frac{1}{2}\left(\frac{\partial F^i}{\partial q_j}-
\frac{\partial F^j}{\partial q_i}\right) dq_i \wedge dq_j .  \]
Here, we use the skew-symmetry of $dq_i \wedge dq_j$ and $dp_i \wedge dp_j$
with respect to the index $i$ and $j$.
It is obvious that the conditions 
(\ref{HC1}) - (\ref{HC3}) lead to the equation
$d\beta=0$. \\

Some of Hamiltonian systems can be defined by the unique function. \\

{\bf Proposition 5}. 
{\it A dynamical system (\ref{eq1}) on the phase space $R^{2n}$ 
is a Hamiltonian system that is defined by Hamiltonian $H=H(q,p)$
if the differential 1-form 
\[ \beta=G^idp_i-F^idq_i, \]
is an exact form $\beta=dH$, 
where $d$ is the exterior derivative and 
$H=H(q,p)$ is a continuous differentiable unique 
function on the phase space. } \\

{\bf Proof.} 
Suppose that the differential 1-form $\beta$, 
which is defined by equation (\ref{beta}), has the form
\[ \beta=dH=\frac{\partial H}{\partial p_i} dp_i+
\frac{\partial H}{\partial q_i} dq_i. \]
In this case, vector fields $(G^i,F^i)$ can be represented
in the form
\[ G^i(q,p)=\frac{\partial H}{\partial p_i}, \quad
F^i(q,p)=-\frac{\partial H}{\partial q_i} .\]
If $H=H(q,p)$ is a continuous differentiable function,
then condition (\ref{HC1}) - (\ref{HC3})
are satisfied. Using proposition 4, we get that this system
is a Hamiltonian system.
The equations of motion for the Hamiltonian system 
(\ref{eq1}) can be written in the form
\be \label{eq2} \frac{dq_i}{dt}=\frac{\partial H}{\partial p_i}, 
\quad \frac{dp_i}{dt}=-\frac{\partial H}{\partial q_i} , \ee
which is uniquely defined by the Hamiltonian $H$.\\

If the exact differential 1-form $\beta$ is equal to zero
($dH=0$), then the equation 
\be \label{HC} H(q,p)-C=0 \ee
defines the stationary states of the Hamiltonian system (\ref{eq1}).
Here, $C$ is a constant.

\section{Fractional Hamiltonian systems}

Fractional generalization of the differential form (\ref{beta}), 
which is used in definition of the Hamiltonian system, 
can be defined in the following form: 
\[ \beta_{\alpha}=G^i (dp_i)^{\alpha}-F^i (dq_i)^{\alpha}. \]
Let us consider the canonical coordinates 
$(x^1,...,x^n,x^{n+1},...,x^{2n})=(q^1,...,q^n,p^{1},...,p^{n})$
in the phase space $R^{2n}$ and
a dynamical system that is defined by the equations 
\be \label{feq1} \frac{dq_i}{dt}=G^i(q,p), 
\quad \frac{dp_i}{dt}=F^i(q,p) . \ee
The fractional generalization of Hamiltonian systems can be defined  
by using fractional generalization of differential forms \cite{FDF}. \\

{\bf Definition 4.} 
{\it A dynamical system (\ref{feq1}) on the phase space $R^{2n}$ 
is called a fractional Hamiltonian system if the 
fractional differential 1-form 
\[ \beta_{\alpha}=G^i (dp_i)^{\alpha}-F^i (dq_i)^{\alpha}, \]
is a closed fractional form 
\be \label{ba} d^{\alpha} \beta_{\alpha}=0 , \ee
where $d^{\alpha}$ is the fractional exterior derivative. 
A dynamical system is called a fractional non-Hamiltonian system 
if the fractional differential 1-form $\beta_{\alpha}$ 
is a nonclosed fractional form $d^{\alpha} \beta_{\alpha} \not=0$.} \\

The fractional exterior derivative for the phase space $R^{2n}$
is defined as 
\be \label{fed2}  d^{\alpha}=(dq_i)^{\alpha} 
{\bf D}^{\alpha}_{q_i}+
(dp_i)^{\alpha} {\bf D}^{\alpha}_{p_i} . \ee
For example, the fractional exterior derivative 
of order $\alpha$ of $q^k$, with the initial point taken 
to be zero and $n=2$, is given by
\be \label{daq2} d^{\alpha} q^k= (dq)^{\alpha} {\bf D}^{\alpha}_{q} q^k+
(dp)^{\alpha} {\bf D}^{\alpha}_{p} q^k . \ee
Using  equations (\ref{xk}) and (\ref{const}), 
we have the following relation 
for the fractional exterior derivative (\ref{fed2}): 
\[ d^{\alpha} q^k= (dq)^{\alpha} 
\frac{\Gamma(k+1) q^{k-\alpha}}{\Gamma(k+1-\alpha)} +
(dp)^{\alpha} \frac{q^{k} p^{-\alpha}}{\Gamma(1-\alpha)} . \]
Let us consider a fractional generalization of the Helmholtz conditions. \\

{\bf Proposition 6.} 
{\it If the right-hand sides of equations (\ref{feq1})
satisfy the fractional generalization of the Helmholtz conditions 
in the following form 
\be \label{HC1f} {\bf D}^{\alpha}_{p_j} G^i-
{\bf D}^{\alpha}_{p_i} G^j= 0, \ee
\be \label{HC2f} {\bf D}^{\alpha}_{q_i} G^j+
{\bf D}^{\alpha}_{p_j} F^i=0, \ee
\be \label{HC3f} {\bf D}^{\alpha}_{q_j} F^i-
{\bf D}^{\alpha}_{q_i} F^j= 0, \ee
then dynamical system (\ref{feq1}) is a 
fractional Hamiltonian system.} \\

{\bf Proof.} 
In the canonical coordinates $(q,p)$, the vector fields that 
define the system have the components $(G^i,F^i)$, 
which are used in equation (\ref{eq1}). 
The 1-form $\beta_{\alpha}$ is defined by the equation
\be \label{beta-a} 
\beta_{\alpha}=G^i(dp_i)^{\alpha}-F^i (dq_i)^{\alpha}. \ee
The exterior derivative for this form can now be given by the relation
\[ d^{\alpha}\beta_{\alpha}=d^{\alpha}(G^i(dp_i)^{\alpha})-
d^{\alpha}(F^i(dq_i)^{\alpha}) .  \]
Using the rule 
\[ {\bf D}^{\alpha}_{x} (fg) =\sum^{\infty}_{k=0} 
\left(^{\alpha}_k \right) ({\bf D}^{\alpha-k}_x f )
\frac{\partial^k g}{\partial x^k} ,\]
and the relation
\[  \frac{\partial^k}{\partial x^k} 
\left((dx)^{\alpha}\right)=0 \quad (k\ge 1) , \]
we get that
\[ d^{\alpha} (A^i (dx_i)^{\alpha})=\sum^{\infty}_{k=0}
(dx_j)^{\alpha} \wedge 
\left(^{\alpha}_k \right)
({\bf D}^{\alpha-k}_{x_j} A^i) 
\frac{\partial^k}{\partial x^k_j} (dx_i)^{\alpha}
=(dx_j)^{\alpha} \wedge (dx_i)^{\alpha}
\left(^{\alpha}_0 \right)
({\bf D}^{\alpha}_{x_j} A^i).  \]
Here, we use
\[ \left(^{\alpha}_k \right)
=\frac{(-1)^{k-1} \alpha \Gamma(k-\alpha)}{\Gamma(1-\alpha) \Gamma(k+1)} . \]
Therefore, we have
\[ d^{\alpha} \beta_{\alpha}= 
{\bf D}^{\alpha}_{q_j} G^i (dq_j)^{\alpha} \wedge (dp_i)^{\alpha}
+{\bf D}^{\alpha}_{p_j} G^i (dp_j)^{\alpha} \wedge (dp_i)^{\alpha}-\]
\[ - {\bf D}^{\alpha}_{q_j} F^i (dq_j)^{\alpha} \wedge (dq_i)^{\alpha}-
{\bf D}^{\alpha}_{p_j} F^i (dp_j)^{\alpha} \wedge (dq_i)^{\alpha} . \]
This equation can be rewritten in an equivalent form
\[ d^{\alpha}\beta_{\alpha}=
\Bigl( {\bf D}^{\alpha}_{q_i} G^j 
+{\bf D}^{\alpha}_{p_j} F^i \Bigr) (dq_i)^{\alpha} \wedge (dp_j)^{\alpha}+
\frac{1}{2}\Bigl( {\bf D}^{\alpha}_{p_i} G^j 
-{\bf D}^{\alpha}_{p_j} G^i \Bigr) (dp_i)^{\alpha} \wedge (dp_j)^{\alpha}+\]
\[ +\frac{1}{2}\Bigl( {\bf D}^{\alpha}_{q_j} F^i-
{\bf D}^{\alpha}_{q_i} F^j \Bigr) (dq_i)^{\alpha} \wedge (dq_j)^{\alpha} .  \]
Here, we use the skew symmetry of $\wedge$.
It is obvious that conditions (\ref{HC1f}) - (\ref{HC3f}) 
lead to the equation $d^{\alpha}\beta_{\alpha}=0$,i.e.,
$\beta_{\alpha}$ is a closed fractional form. \\


Let us define the Hamiltonian for the fractional Hamiltonian systems.\\

{\bf Proposition 7.} 
{\it A dynamical system (\ref{feq1}) on the phase space  $R^{2n}$ 
is a fractional Hamiltonian system that is defined by the Hamiltonian 
$H=H(q,p)$ if the fractional differential 1-form 
\[ \beta_{\alpha}=G^i (dp_i)^{\alpha}-F^i (dq_i)^{\alpha}, \]
is an exact fractional form
\be \label{b2a}  \beta_{\alpha}=d^{\alpha} H , \ee
where $d^{\alpha}$ is the fractional exterior derivative 
and $H=H(q,p)$ is a continuous differentiable 
function on the phase space. }\\

{\bf Proof}. 
Suppose that the fractional differential 1-form $\beta_{\alpha}$, 
which is defined by equation (\ref{beta-a}), has the form
\[ \beta_{\alpha}=d^{\alpha}H=(dp_i)^{\alpha}
{\bf D}^{\alpha}_{p_i} H +(dq_i)^{\alpha} {\bf D}^{\alpha}_{q_i} H . \]
In this case, vector fields $(G^i,F^i)$ can be represented
in the form
\[ G^i(q,p)={\bf D}^{\alpha}_{p_i} H, \quad
F^i(q,p)=-{\bf D}^{\alpha}_{q_i} H .\]
Therefore, the equations of motion for fractional Hamiltonian systems 
can be written in the form
\be \label{feq2} \frac{dq_i}{dt}={\bf D}^{\alpha}_{p_i} H, 
\quad \frac{dp_i}{dt}=-{\bf D}^{\alpha}_{q_i} H . \ee

The fractional differential 1-form $\beta_{\alpha}$ 
for the fractional Hamiltonian system with Hamiltonian $H$ 
can be written in the form $\beta_{\alpha}=d^{\alpha}H$. 
If the exact fractional differential 1-form $\beta_{\alpha}$ 
is equal to zero ($d^{\alpha}H=0$), then we can get the equation 
that defines the stationary states of the Hamiltonian system. \\

{\bf Proposition 8.}
{\it The stationary states of the fractional Hamiltonian system (\ref{feq2})
are defined by the equation
\be \label{ssHs} H(q,p)-|\prod^n_{i=1} q_i p_i|^{\alpha-m} 
\sum^{m-1}_{k_1=0,l_1=0} ... \sum^{m-1}_{k_n=0,l_n=0} 
C_{k_1...k_n l_1...l_n} \prod^n_{i=1}(q_i)^{k_i} (p_i)^{l_i}=0 , \ee
where $C_{k_1...k_n,l_1,...,l_n}$ are constants and $m$ is the first whole 
number greater than  or equal to $\alpha$.} \\
 
{\bf Proof}. This proposition is a corollary of proposition 3. \\

\section{Example of fractional Hamiltonian system}

Let us consider a dynamical system in phase space $R^2$ ($n=1$)
that is defined by the equation
\be \label{fo} \frac{dq}{dt}={\bf D}^{\alpha}_{p} H, 
\quad \frac{dp}{dt}=-{\bf D}^{\alpha}_{q} H , \ee
where the fractional order $0< \alpha\le 1$ and
the Hamiltonian $H(q,p)$ has the form
\be \label{Hqp} H(q,p)=ap^2+bq^2 . \ee
If $\alpha=1$, then equation (\ref{fo}) describes the 
linear harmonic oscillator. 

If the exact fractional differential 1-form 
\[ \beta_{\alpha}=d^{\alpha}H=(dp)^{\alpha}
{\bf D}^{\alpha}_{p} H +(dq)^{\alpha} {\bf D}^{\alpha}_{q} H  \] 
is equal to zero ($d^{\alpha}H=0$), then the equation 
\[ H(q,p)-C |qp|^{\alpha-1}=0 \]
defines the stationary states of the system (\ref{fo}). 
Here, $C$ is a constant. If $\alpha=1$, we get  the usual
stationary-state equation (\ref{HC}).
 
Using equation (\ref{Hqp}), we get the following equation 
for stationary states:
\be \label{fabc} |qp|^{1-\alpha} (ap^2+bq^2)=C . \ee
If $\alpha=1$, then we get the equation $ap^2+bq^2=C$, 
which describes the ellipse.

\section{Conclusion}

Fractional derivatives and integrals \cite{SKM,OS} have 
found many applications in recent studies in physics.
The interest in fractional analysis has been growing continually 
during the past few years [3-17].
Using the fractional derivatives and fractional differential
forms, we consider the fractional generalization
of gradient and Hamiltonian systems.
In the general case, the fractional gradient and Hamiltonian systems cannot
be considered as gradient and Hamiltonian systems. 
The class of fractional gradient and Hamiltonian systems
is wider than the usual class of gradient and Hamiltonian 
dynamical systems. 
The gradient and Hamiltonian systems can be considered as
aspecial case of fractional gradient and Hamiltonian systems. 
Therefore, it is possible to generalize the application 
of catastrophe and bifurcation theory from
gradient to a wider class of fractional gradient dynamical systems. 
Note that quantization of the fractional Hamiltonian systems
can be realized by the method suggested 
in \cite{Tarpla1,Tarmsu,Tarsam,JPA04}.






\end{document}